\documentclass[11pt]{article}
\usepackage[a4paper,margin=2.45cm]{geometry}
\usepackage{iftex}
\usepackage[expansion=false]{microtype}
\ifPDFTeX
  \usepackage[T1]{fontenc}
  \usepackage[utf8]{inputenc}
\else
  \usepackage{fontspec}
\fi
\usepackage[english]{babel}
\usepackage{amsmath,amssymb,amsfonts,amsthm,mathtools}
\ifPDFTeX\else
  \usepackage{unicode-math}
\fi
\usepackage{enumitem}
\usepackage{microtype}
\usepackage{xcolor}
\usepackage{hyperref}
\usepackage[nameinlink,capitalize]{cleveref}
\hypersetup{
 colorlinks=true,
 linkcolor=blue!45!black,
 citecolor=blue!45!black,
 urlcolor=blue!45!black,
 pdftitle={Forward-Backward Stochastic Differential Equations with Time-Delayed Generators under an Osgood Condition},
 pdfauthor={Harouna Coulibaly}
}
\newtheorem{theorem}{Theorem}[section]
\newtheorem{proposition}[theorem]{Proposition}
\newtheorem{lemma}[theorem]{Lemma}
\newtheorem{corollary}[theorem]{Corollary}
\newtheorem{definition}[theorem]{Definition}
\newtheorem{example}[theorem]{Example}
\theoremstyle{remark}
\newtheorem{remark}[theorem]{Remark}

\newcommand{\E}{\mathbb E}
\newcommand{\Pp}{\mathbb P}
\newcommand{\R}{\mathbb R}
\newcommand{\F}{\mathbb F}
\newcommand{\Sp}{\mathcal S^2}
\newcommand{\Hp}{\mathcal H^2}
\newcommand{\norm}[1]{\left\lVert#1\right\rVert}

\makeatletter
\renewcommand{\@biblabel}[1]{#1.}
\makeatother

\title{Forward-Backward Stochastic Differential Equations\\
with Time-Delayed Generators under an Osgood Condition}
\author{
Harouna Coulibaly\\
\small École Supérieure Africaine des Technologies de l'Information\\
\small et de la Communication (ESATIC),\\
\small 18 BP 2202 Abidjan, Côte d'Ivoire\\
\small \texttt{harounbolz@yahoo.fr}
}
\date{}

\begin{document}
\maketitle

\begin{abstract}
We study a class of forward-backward stochastic differential equations in
which the forward coefficients and the backward driver depend on delayed
segments $X_t$ of the forward state, while $Y(t)$ and $Z(t)$ remain current
arguments of the driver. The forward coefficients satisfy an
integral Osgood condition on the path space, while the backward driver
satisfies a one-sided Osgood condition in the backward state and a
Lipschitz condition in the martingale integrand. The system is decoupled in
the forward direction, but the backward equation retains the full delayed
forward history. We establish delay estimates, a Bihari--Osgood lemma,
existence by successive approximation for the forward equation, pathwise
uniqueness, and quantitative stability. Backward existence is obtained from
the theory of weakly monotone BSDEs after verifying its local integrability
requirements. The assumptions contain the globally Lipschitz case, but also
allow weakly monotone drivers and moduli such as
$r\mapsto r\log(1/r)$ near the origin. The distinction between this
decoupled Osgood framework and the fully coupled Lipschitz theory is made
explicit.


\end{abstract}

\section{Introduction}

Forward-backward stochastic differential equations (FBSDEs) arise in
stochastic control, mathematical finance and probabilistic representations
of nonlinear partial differential equations. The classical well-posedness
theory for backward stochastic differential equations was initiated by
Pardoux and Peng \cite{PardouxPeng1990}. For coupled FBSDEs, contraction
arguments on short intervals and monotonicity methods provide two standard
approaches; see Antonelli \cite{Antonelli1993}, Hu and Peng
\cite{HuPeng1995}, and Peng and Wu \cite{PengWu1999}.

Time-delayed generators describe systems whose coefficients depend on past
states. Delong and Imkeller \cite{DelongImkeller2010} established a basic
theory for delayed BSDEs. Aman, Coulibaly and Djordjević
\cite{AmanCoulibalyDjordjevic2023} treated coupled FBSDEs with delayed
generators under a delayed Lipschitz condition and an explicit smallness
restriction. Coulibaly and Aman \cite{CoulibalyAman2020} also
studied non-Lipschitz BSDEs in which past backward values enter the
generator. Here, by contrast, $Y(t)$ and $Z(t)$ are current values and the
backward delay acts only through the forward segment $X_t$.

The Osgood condition is a natural weakening of global Lipschitz continuity.
Its decisive feature is the divergence
\[
\int_{0^+}\frac{dr}{\rho(r)}=+\infty,
\]
which permits uniqueness through Bihari's inequality. Mao
\cite{Mao1995} used this mechanism for BSDEs with non-Lipschitz
coefficients. Fan, Jiang and Davison \cite{FanJiangDavison2013} developed
an Osgood theory in which the driver is Osgood in the backward variable and
Lipschitz in the martingale integrand.
For weakly monotone BSDEs with general local growth in the backward state,
we use Fan's $L^p$ theory \cite{Fan2015}. The forward equation is a
stochastic functional differential equation with bounded memory; related
results include von Renesse and Scheutzow \cite{VonRenesseScheutzow2010},
while Picard--Bihari arguments for functional equations also appear in Wei
and Cai \cite{WeiCai2013}.

The objective of this paper is to combine an Osgood modulus with a delayed
path dependence. We consider
\begin{align}
X(t)&=\varphi(0)+\int_0^t b(s,X_s)\,ds
       +\int_0^t\sigma(s,X_s)\,dW(s), \label{eq:forward}\\
Y(t)&=\xi+\int_t^T f(s,X_s,Y(s),Z(s))\,ds
       -\int_t^T Z(s)\,dW(s), \label{eq:backward}
\end{align}
where $X_s=(X(s+u))_{-\delta\le u\le0}$. The forward coefficients depend
on the whole delayed segment. The backward driver is also time-delayed
through $X_s$, is Osgood in $Y$, and is Lipschitz in $Z$.

The forward equation is deliberately decoupled from $(Y,Z)$. This is not a
cosmetic restriction. In the fully coupled system of
\cite{AmanCoulibalyDjordjevic2023}, a Lipschitz contraction closes the
simultaneous estimates of $(X,Y,Z)$. The present Picard--Bihari argument
treats the triangular system by solving the two equations successively; a
fully coupled Osgood extension would require a separate argument. The
present formulation isolates this rigorous non-Lipschitz contribution.

Our main contributions are:
\begin{enumerate}[label=(\roman*)]
\item an Osgood delay estimate on the segment space;
\item existence and uniqueness of the delayed forward equation;
\item existence and uniqueness of the backward equation under a one-sided
Osgood condition in $Y$ and a Lipschitz condition in $Z$;
\item stability with respect to the initial path, terminal condition and
coefficients;
\item quantitative continuous dependence obtained by composing the forward
and backward Osgood controls.
\end{enumerate}

\section{Probability space and functional setting}

Let $(\Omega,\mathcal F,\F)=((\mathcal F_t)_{0\le t\le T},\Pp)$ be a
complete filtered probability space carrying a $d$-dimensional Brownian
motion $W$. The filtration is the usual augmentation of the natural
filtration of $W$.

For a Euclidean space $E$, let $\Sp(E)$ be the space of continuous adapted
$E$-valued processes $U$ such that
\[
\norm{U}_{\Sp(E)}^2:=\E\sup_{0\le t\le T}|U(t)|^2<\infty.
\]
Let $\Hp(\R^{m\times d})$ be the space of progressively measurable
processes $Z$ satisfying
\[
\norm{Z}_{\Hp}^2:=\E\int_0^T|Z(t)|^2dt<\infty.
\]

Fix a delay length $\delta\in(0,T]$ and a probability measure $\alpha$ on
$[-\delta,0]$. An initial forward path
$\varphi\in L^2(\Omega;C([-\delta,0];\R^n))$ is assumed
$\mathcal F_0$-measurable. For any continuous process $X$ on
$[-\delta,T]$, define
\[
X_t(u):=X(t+u),\qquad -\delta\le u\le0,
\]
and
\[
\norm{X_t-X'_t}_{\alpha}^2
:=\int_{-\delta}^0|X(t+u)-X'(t+u)|^2\alpha(du).
\]
On $C([-\delta,0];\R^n)$, the quantity $\norm{\eta}_\alpha$ is in
general a seminorm; it is a norm on the quotient space $L^2(\alpha)$, and
also on the continuous segment space when $\alpha$ has full support. No
full-support assumption is needed below, since the proofs use only
\[
\norm{\eta}_\alpha\le\norm{\eta}_\infty.
\]
We use the prescribed history
\begin{equation}
X(t)=\varphi(t),\qquad -\delta\le t\le0. \label{eq:history}
\end{equation}

\begin{definition}[Osgood modulus]\label{def:osgood}
A function $\rho:[0,\infty)\to[0,\infty)$ is called an Osgood modulus if
it is continuous, nondecreasing and concave, $\rho(0)=0$,
$\rho(r)>0$ for $r>0$, and
\begin{equation}
\int_{0^+}\frac{dr}{\rho(r)}=+\infty. \label{eq:osgood}
\end{equation}
\end{definition}

Concavity and $\rho(0)=0$ imply the linear-growth estimate
\begin{equation}
\rho(r)\le c_\rho(1+r),\qquad r\ge0, \label{eq:rho-growth}
\end{equation}
for a constant $c_\rho>0$.

\subsection{Assumptions}

We impose the following conditions.
\begin{description}[leftmargin=1.35cm,labelsep=.55em]
\item[(O1)] The initial path satisfies
$\E\norm{\varphi}_{\infty}^2<\infty$ and
$\xi\in L^2(\Omega,\mathcal F_T;\R^m)$.

\item[(O2)] For every continuous adapted input $X$, the processes
$b(t,X_t)$ and $\sigma(t,X_t)$ are progressively measurable. Moreover,
\[
\E\int_0^T\bigl(|b(t,0)|^2+|\sigma(t,0)|^2\bigr)dt<\infty.
\]
There is an Osgood modulus $\rho_X$ and a constant $L_X>0$ such that
\begin{equation}
|b(t,\eta)-b(t,\eta')|^2+
|\sigma(t,\eta)-\sigma(t,\eta')|^2
\le L_X\rho_X\!\left(\norm{\eta-\eta'}_\alpha^2\right)
\label{eq:forward-osgood}
\end{equation}
for all admissible segments, $dt\otimes d\Pp$-a.e.

\item[(O3)] For every triple
\[
(X,Y,Z)\in\Sp(\R^n)\times\Sp(\R^m)\times\Hp(\R^{m\times d}),
\]
the process
$f(t,X_t,Y(t),Z(t))$ is progressively measurable and
\[
\E\int_0^T|f(t,0,0,0)|^2dt<\infty.
\]
There are an Osgood modulus $\rho_Y$ and constants $L_f,L_z>0$ such that
\begin{align}
|f(t,\eta,y,z)-f(t,\eta',y,z)|^2
&\le L_f\rho_X\!\left(\norm{\eta-\eta'}_\alpha^2\right),\label{eq:f-x}\\
2\langle y-y',f(t,\eta,y,z)-f(t,\eta,y',z)\rangle
&\le L_f\rho_Y(|y-y'|^2), \label{eq:f-y}\\
|f(t,\eta,y,z)-f(t,\eta,y,z')|
&\le L_z|z-z'|. \label{eq:f-z}
\end{align}

\item[(O4)] For almost every $(\omega,t)$, the maps
$\eta\mapsto b(t,\eta)$, $\eta\mapsto\sigma(t,\eta)$ and
$(\eta,y,z)\mapsto f(t,\eta,y,z)$ are continuous.

\item[(O5)] The driver has locally integrable growth in the backward state:
for every $R>0$,
\begin{equation}
\E\int_0^T\sup_{|y|\le R}
|f(t,0,y,0)-f(t,0,0,0)|\,dt<\infty .
\label{eq:local-growth}
\end{equation}
\end{description}

Condition \emph{(O5)} supplies a control that is essential for a theorem on
the whole class considered here: the one-sided inequality
\eqref{eq:f-y} controls an energy product but not the norm of the driver.
For example, $f(t,y)=-y/(T-t)$ for $t<T$, with $f(T,y)=0$, satisfies
\eqref{eq:f-y} with a nonpositive left-hand side, but prevents a continuous
solution with a nonzero terminal value because its drift is not locally
integrable near $T$. This shows the need for an additional growth condition
in the general theorem; it does not assert that \emph{(O5)} is necessary for
every individual datum (for instance, the same driver with $\xi=0$ admits
the zero solution).

\begin{remark}
Taking $\rho_X(r)=\rho_Y(r)=r$ makes the segment condition Lipschitz and
the condition in $y$ one-sided Lipschitz; it does not by itself impose a
two-sided global Lipschitz bound in $y$. For example, $f(y)=-y^3$ satisfies
\eqref{eq:f-y} with this linear modulus and satisfies \emph{(O5)}, although
it is not globally Lipschitz. A genuinely nonlinear Osgood modulus is
obtained from
\[
\rho(r)=
\begin{cases}
r\log(1/r),&0<r\le e^{-1},\\
e^{-1},&r>e^{-1}.
\end{cases}
\]
Indeed, $\int_{0^+}dr/[r\log(1/r)]=+\infty$.
\end{remark}

\section{Preliminary estimates}

\begin{lemma}[Delay estimate]\label{lem:delay}
Let $X,X'$ be continuous processes on $[-\delta,T]$. Then, for every
$t\in[0,T]$,
\begin{equation}
\int_0^t\E\norm{X_s-X'_s}_{\alpha}^2ds
\le\int_0^t\E\sup_{-\delta\le r\le s}|X(r)-X'(r)|^2ds.
\label{eq:delay-estimate}
\end{equation}
Moreover, for every concave Osgood modulus $\rho$,
\begin{equation}
\int_0^t\E\rho\!\left(\norm{X_s-X'_s}_{\alpha}^2\right)ds
\le\int_0^t\rho\!\left(
\E\sup_{-\delta\le r\le s}|X(r)-X'(r)|^2\right)ds.
\label{eq:delay-osgood}
\end{equation}
\end{lemma}

\begin{proof}
For $u\in[-\delta,0]$, one has $-\delta\le s+u\le s$, and hence
\[
|X(s+u)-X'(s+u)|^2
\le\sup_{-\delta\le r\le s}|X(r)-X'(r)|^2.
\]
Integration with respect to the probability measure $\alpha$ gives
\[
\norm{X_s-X'_s}_{\alpha}^2
\le\sup_{-\delta\le r\le s}|X(r)-X'(r)|^2.
\]
Taking expectations and integrating in $s$ proves
\eqref{eq:delay-estimate}. Since $\rho$ is nondecreasing,
\[
\rho(\norm{X_s-X'_s}_{\alpha}^2)
\le\rho\!\left(\sup_{-\delta\le r\le s}|X(r)-X'(r)|^2\right).
\]
Jensen's inequality for the concave function $\rho$ yields
\[
\E\rho\!\left(\sup_{-\delta\le r\le s}|X(r)-X'(r)|^2\right)
\le\rho\!\left(\E\sup_{-\delta\le r\le s}|X(r)-X'(r)|^2\right).
\]
Integration proves \eqref{eq:delay-osgood}.
\end{proof}

\begin{lemma}[Bihari--Osgood]\label{lem:bihari}
Let $\rho$ be an Osgood modulus and let $u:[0,T]\to[0,\infty)$ be bounded
and measurable. If
\begin{equation}
u(t)\le a+C\int_0^t\rho(u(s))ds,\qquad 0\le t\le T, \label{eq:bihari-forward}
\end{equation}
then, if $a>0$, for every $r_0>0$,
\begin{equation}
G_{r_0}(u(t))\le G_{r_0}(a)+Ct,\qquad
G_{r_0}(r):=\int_{r_0}^{r}\frac{dq}{\rho(q)}, \label{eq:bihari-bound}
\end{equation}
where $G_{r_0}$ is extended by $G_{r_0}(0)=-\infty$,
and $u(t)\le G_{r_0}^{-1}(G_{r_0}(a)+Ct)$. If $a=0$, then
$u\equiv0$. The same conclusions hold for the
backward inequality
\[
u(t)\le a+C\int_t^T\rho(u(s))ds.
\]
\end{lemma}

\begin{proof}
For $\varepsilon>0$, put
\[
v_\varepsilon(t):=a+\varepsilon+
C\int_0^t\rho(u(s))ds.
\]
Then $u(t)\le v_\varepsilon(t)$ and, because $\rho$ is nondecreasing,
\[
v_\varepsilon'(t)=C\rho(u(t))
\le C\rho(v_\varepsilon(t))
\quad\text{for a.e. }t.
\]
Consequently,
\[
\frac{d}{dt}G_{r_0}(v_\varepsilon(t))
=\frac{v_\varepsilon'(t)}{\rho(v_\varepsilon(t))}
\le C.
\]
Integration from $0$ to $t$ gives
\[
G_{r_0}(v_\varepsilon(t))
\le G_{r_0}(a+\varepsilon)+Ct.
\]
Letting $\varepsilon\downarrow0$ proves \eqref{eq:bihari-bound} when
$a>0$. If $a=0$ and $u(t_0)>0$ for some $t_0$, repeat the argument after
adding $\varepsilon$ and let $\varepsilon\downarrow0$. The left-hand
integral then crosses every interval $(\varepsilon,u(t_0))$, contradicting
\eqref{eq:osgood}; hence $u=0$. Moreover, \eqref{eq:rho-growth} implies
$G_{r_0}(r)\to+\infty$ as $r\to+\infty$, so the inverse in the statement
is defined on the required range. For the backward form, apply the forward
argument to $\widetilde u(r)=u(T-r)$.
\end{proof}

\begin{lemma}[Basic BSDE estimate]\label{lem:bsde}
Let
\[
Y(t)=\eta+\int_t^T h(s)ds-\int_t^TZ(s)dW(s),
\]
where $\eta\in L^2(\mathcal F_T)$, $h$ is progressively measurable, and
$\E\int_0^T|h(s)|^2ds<\infty$. Then
\begin{equation}
\norm{Y}_{\Sp}^2+\norm{Z}_{\Hp}^2
\le C_T\E\left[|\eta|^2+\int_0^T|h(s)|^2ds\right].
\label{eq:linear-bsde}
\end{equation}
The same estimate holds for differences.
\end{lemma}

\begin{proof}
Set $\zeta=\eta+\int_0^Th(s)ds$ and
$M(t)=\E[\zeta\mid\mathcal F_t]$. The martingale representation theorem
gives
\[
M(t)=M(0)+\int_0^tZ(s)dW(s).
\]
Define $Y(t)=M(t)-\int_0^th(s)ds$. This is the unique solution. Doob's
inequality, Cauchy--Schwarz and the BDG inequalities give
\begin{align*}
\E\sup_{t\le T}|Y(t)|^2
&\le C\E|\zeta|^2+C\E\int_0^T|h(s)|^2ds,\\
\E\int_0^T|Z(s)|^2ds
&\le C\E|\zeta|^2.
\end{align*}
Finally,
\[
\E|\zeta|^2\le2\E|\eta|^2+2T\E\int_0^T|h(s)|^2ds.
\]
Combining these bounds proves \eqref{eq:linear-bsde}. Subtraction proves
the difference estimate.
\end{proof}

\section{The delayed forward equation}

\begin{proposition}[A priori estimate]\label{prop:forward-apriori}
Under \emph{(O1)--(O2)}, every solution of \eqref{eq:forward} satisfies
\begin{equation}
\E\sup_{-\delta\le t\le T}|X(t)|^2
\le C\left(1+\E\norm{\varphi}_\infty^2+
\E\int_0^T(|b(t,0)|^2+|\sigma(t,0)|^2)dt\right),
\label{eq:forward-apriori}
\end{equation}
where $C$ depends only on $T,L_X$ and $\rho_X$ through the linear-growth
constant in \eqref{eq:rho-growth}.
\end{proposition}

\begin{proof}
For $t\in[0,T]$, the elementary inequality, Cauchy--Schwarz and BDG give
\begin{align*}
\E\sup_{0\le r\le t}|X(r)|^2
&\le3\E|\varphi(0)|^2
+3t\E\int_0^t|b(s,X_s)|^2ds\\
&\quad+C\E\int_0^t|\sigma(s,X_s)|^2ds.
\end{align*}
By \emph{(O2)}, \eqref{eq:rho-growth}, and
\[
\norm{X_s}_{\alpha}^2
\le\norm{\varphi}_{\infty}^2+
\sup_{0\le r\le s}|X(r)|^2,
\]
we have
\begin{align*}
|b(s,X_s)|^2+|\sigma(s,X_s)|^2
&\le2|b(s,0)|^2+2|\sigma(s,0)|^2\\
&\quad+2L_Xc_{\rho_X}
\left(1+\norm{\varphi}_\infty^2+
\sup_{0\le r\le s}|X(r)|^2\right).
\end{align*}
Substitution and Gronwall's lemma prove the bound on $[0,T]$; adding the
prescribed history proves \eqref{eq:forward-apriori}.
\end{proof}

\begin{theorem}[Forward Osgood well-posedness]\label{thm:forward}
Under \emph{(O1), (O2), (O4)}, equation \eqref{eq:forward} has a unique
solution $X\in\Sp(\R^n)$, extended by \eqref{eq:history}.
\end{theorem}

\begin{proof}
\emph{Uniqueness.}
Let $X,X'$ be two solutions with the same history, and put
\[
u(t):=\E\sup_{0\le r\le t}|X(r)-X'(r)|^2.
\]
BDG, Cauchy--Schwarz and \eqref{eq:forward-osgood} yield
\[
u(t)\le C_{T,L_X}\int_0^t
\E\rho_X(\norm{X_s-X'_s}_{\alpha}^2)ds.
\]
By Lemma~\ref{lem:delay},
\[
u(t)\le C_{T,L_X}\int_0^t\rho_X(u(s))ds.
\]
Lemma~\ref{lem:bihari} gives $u(t)=0$ for all $t$, proving uniqueness.

\emph{Existence.}
Set $X^0(t)=\varphi(0)$ for $t\ge0$ and $X^0(t)=\varphi(t)$ for
$t\in[-\delta,0]$. Recursively define, for $j\ge0$,
\[
X^{j+1}(t)=\varphi(0)+\int_0^t b(s,X^j_s)ds
 +\int_0^t\sigma(s,X^j_s)dW(s),\qquad t\ge0,
\]
and prescribe the same history $\varphi$ on $[-\delta,0]$. Progressive
measurability in \emph{(O2)} and the linear growth following from
\eqref{eq:forward-osgood} and \eqref{eq:rho-growth} show inductively that
every iterate is well defined and belongs to $\Sp$.

Put
\[
A_j(t):=\E\sup_{-\delta\le r\le t}|X^j(r)|^2.
\]
The BDG and Cauchy--Schwarz estimates used in Proposition
\ref{prop:forward-apriori} give constants $C_0,C_1$, independent of $j$,
such that
\[
A_{j+1}(t)\le C_0+C_1\int_0^t A_j(s)ds,
\qquad C_0\ge\E\norm{\varphi}_\infty^2.
\]
An induction on $j$, using the exponential series, yields
\[
\sup_{j\ge0}\sup_{t\le T}A_j(t)\le C_0e^{C_1T}.
\]

For $N\ge0$, define
\[
U_N(t):=\sup_{j,k\ge N}\E\sup_{0\le r\le t}
|X^j(r)-X^k(r)|^2.
\]
Since all iterates have the same negative-time history, BDG,
\eqref{eq:forward-osgood}, Jensen's inequality and Lemma~\ref{lem:delay}
give
\[
U_{N+1}(t)\le C\int_0^t\rho_X(U_N(s))ds.
\]
The functions $U_N$ decrease in $N$ and are uniformly bounded by the
preceding moment estimate. Hence $U_N(t)\downarrow U(t)$ and dominated
convergence, together with continuity of $\rho_X$, gives
\[
U(t)\le C\int_0^t\rho_X(U(s))ds.
\]
Lemma~\ref{lem:bihari} implies $U\equiv0$. Thus $(X^j)$ is Cauchy in
$\Sp$ and converges to some $X$.

Finally, \eqref{eq:forward-osgood}, Jensen and the delay estimate yield
\begin{align*}
\E\int_0^T|b(s,X^j_s)-b(s,X_s)|^2ds
&\le L_XT\rho_X(\norm{X^j-X}_{\Sp}^2),\\
\E\int_0^T|\sigma(s,X^j_s)-\sigma(s,X_s)|^2ds
&\le L_XT\rho_X(\norm{X^j-X}_{\Sp}^2).
\end{align*}
Both right-hand sides tend to zero. Cauchy--Schwarz and BDG therefore
allow passage to the limit in the two integrals defining $X^{j+1}$.
The limit satisfies \eqref{eq:forward}; uniqueness was proved above.
\end{proof}

\section{The backward equation and the FBSDE}

\begin{proposition}[Backward stability estimate]\label{prop:backward-stability}
Let $X,X'\in\Sp(\R^n)$ be continuously extended to $[-\delta,0]$ by
prescribed histories in $L^2(\Omega;C([ -\delta,0];\R^n))$, and let
$(Y,Z)$, $(Y',Z')$ solve the backward equations with terminal values
$\xi,\xi'$ and forward inputs $X,X'$.
Then there is a constant $C=C(T,L_f,L_z)$ such that
\begin{align}
v(t)+\frac12\E\int_t^T|Z(s)-Z'(s)|^2ds
&\le C\E|\xi-\xi'|^2
+C\int_t^T\overline\rho_Y(v(s))ds \notag\\
&\quad+C\int_t^T
\rho_X\!\left(\E\sup_{-\delta\le r\le s}|X(r)-X'(r)|^2\right)ds,
\label{eq:backward-stability}
\end{align}
where $v(t):=\E|Y(t)-Y'(t)|^2$ and
$\overline\rho_Y(r):=r+\rho_Y(r)$.
\end{proposition}

\begin{proof}
Write $\Delta Y=Y-Y'$, $\Delta Z=Z-Z'$ and decompose the driver difference:
\begin{align*}
\Delta f(s)
&=f(s,X_s,Y(s),Z(s))-f(s,X'_s,Y'(s),Z'(s))\\
&=\bigl[f(s,X_s,Y,Z)-f(s,X'_s,Y,Z)\bigr]\\
&\quad+\bigl[f(s,X'_s,Y,Z)-f(s,X'_s,Y',Z)\bigr]\\
&\quad+\bigl[f(s,X'_s,Y',Z)-f(s,X'_s,Y',Z')\bigr].
\end{align*}
Itô's formula, first up to
$\tau_N=\inf\{t:\int_0^t|\Delta Y(s)|^2|\Delta Z(s)|^2ds\ge N\}\wedge T$
and then by localization, gives
\begin{align*}
\E|\Delta Y(t)|^2+\E\int_t^T|\Delta Z(s)|^2ds
&=\E|\xi-\xi'|^2
+2\E\int_t^T\langle\Delta Y(s),\Delta f(s)\rangle ds.
\end{align*}
Indeed, with
\[
M_t=\int_0^t\langle\Delta Y(s),\Delta Z(s)dW(s)\rangle,
\quad A=\E\sup_{t\le T}|\Delta Y(t)|^2,
\quad B=\E\int_0^T|\Delta Z(s)|^2ds,
\]
BDG at exponent one and Cauchy--Schwarz give
$\E\sup_{t\le T}|M_t|\le C\sqrt{AB}<\infty$. Thus the stopped martingales
are uniformly integrable and their expectations vanish in the limit. The
positive part of the drift product is integrable by
\eqref{eq:f-x}--\eqref{eq:f-z} and Young's inequality; the localized Itô
identity then also defines its signed integral.
For the delayed forward term, Young's inequality and \eqref{eq:f-x} give
\[
2|\Delta Y|\,|f(s,X_s,Y,Z)-f(s,X'_s,Y,Z)|
\le|\Delta Y|^2+
L_f\rho_X(\norm{X_s-X'_s}_\alpha^2).
\]
The $Y$-term is controlled by \eqref{eq:f-y}. For the $Z$-term,
\[
2L_z|\Delta Y||\Delta Z|
\le2L_z^2|\Delta Y|^2+\frac12|\Delta Z|^2.
\]
Taking expectations, using concavity and Lemma~\ref{lem:delay}, we obtain
\begin{align*}
v(t)+\frac12\E\int_t^T|\Delta Z(s)|^2ds
&\le\E|\xi-\xi'|^2
+C\int_t^T v(s)ds\\
&\quad+C\int_t^T\rho_Y(v(s))ds\\
&\quad+C\int_t^T
\rho_X\!\left(\E\sup_{-\delta\le r\le s}|X(r)-X'(r)|^2\right)ds.
\end{align*}
The function $\overline\rho_Y(r)=r+\rho_Y(r)$ is continuous, concave and
nondecreasing. Moreover, concavity and positivity imply
$r\le C\rho_Y(r)$ near zero, so
$\overline\rho_Y(r)\le C'\rho_Y(r)$ there. Thus
$\int_{0^+}dr/\overline\rho_Y(r)=\infty$.
Absorbing the linear term into $\overline\rho_Y$ gives
\eqref{eq:backward-stability}.
\end{proof}

\begin{theorem}[Osgood well-posedness of the delayed FBSDE]
\label{thm:main}
Under \emph{(O1)--(O5)}, the system
\eqref{eq:forward}--\eqref{eq:backward} has a unique solution
\[
(X,Y,Z)\in\Sp(\R^n)\times\Sp(\R^m)\times\Hp(\R^{m\times d}).
\]
\end{theorem}

\begin{proof}
Theorem~\ref{thm:forward} gives a unique forward state $X$. It remains to
solve the BSDE with the progressively measurable driver
\[
F_X(t,y,z):=f(t,X_t,y,z).
\]
By \emph{(O3)--(O4)}, this driver is one-sided Osgood in $y$, Lipschitz in
$z$, continuous in $(y,z)$, and
$F_X(\cdot,0,0)\in L^2(dt\otimes d\Pp)$. The last assertion follows from
\begin{align*}
|F_X(t,0,0)|^2
&\le2|f(t,0,0,0)|^2+
2L_f\rho_X(\norm{X_t}_\alpha^2),
\end{align*}
the linear growth of $\rho_X$ and
Proposition~\ref{prop:forward-apriori}.

We now verify the local-growth condition for this fixed forward state. For
each $R>0$, two applications of \eqref{eq:f-x} give
\begin{align*}
\sup_{|y|\le R}|F_X(t,y,0)-F_X(t,0,0)|
&\le \sup_{|y|\le R}|f(t,0,y,0)-f(t,0,0,0)|\\
&\quad+2\sqrt{L_f\rho_X(\norm{X_t}_\alpha^2)}.
\end{align*}
The first term is integrable by \emph{(O5)}. For the second one,
Cauchy--Schwarz, the linear growth of $\rho_X$, and Proposition
\ref{prop:forward-apriori} imply
\[
\E\int_0^T\sqrt{\rho_X(\norm{X_t}_\alpha^2)}dt<\infty.
\]
For fixed $(y,z)$, progressive measurability follows from \emph{(O3)} by
using constant adapted inputs; continuity and approximation by rational
parameters then give the required joint measurability. Thus $F_X$ is
progressively measurable, continuous in $(y,z)$, weakly
monotone in $y$ with the Osgood modulus $\rho_Y$, Lipschitz in $z$, square
integrable at $(0,0)$, and locally integrable in $y$. These are precisely
the hypotheses (H1)--(H5) of Theorem~1 of Fan \cite{Fan2015} for $p=2$:
in Fan's normalization the weak-monotonicity modulus is
$(L_f/2)\rho_Y$, and
\[
\E\left(\int_0^T|F_X(t,0,0)|dt\right)^2
\le T\E\int_0^T|F_X(t,0,0)|^2dt<\infty.
\]
That theorem
therefore yields a pair $(Y,Z)\in\Sp\times\Hp$ satisfying
\eqref{eq:backward}; in particular,
\[
\int_0^T|f(s,X_s,Y(s),Z(s))|ds<\infty\quad\Pp\text{-a.s.}
\]

If $(Y,Z)$ and $(Y',Z')$ are two solutions for the same $X$ and $\xi$,
Proposition~\ref{prop:backward-stability} gives
\[
v(t)\le C\int_t^T\overline\rho_Y(v(s))ds.
\]
The backward Bihari--Osgood lemma gives $v=0$. Returning to Itô's identity
then gives $\norm{Z-Z'}_{\Hp}=0$. This proves uniqueness of the backward
pair and hence of the full triple.
\end{proof}

\section{Stability and continuous dependence}

\begin{theorem}[Stability]\label{thm:stability}
For $i=1,2$, let $(\varphi_i,\xi_i,b_i,\sigma_i,f_i)$ satisfy the preceding
assumptions with common moduli and constants on the same filtered probability
space, driven by the same Brownian motion, with the same horizon $T$, delay
$\delta$, and delay measure $\alpha$. Let $(X_i,Y_i,Z_i)$ be the
corresponding solutions. Define
\begin{align*}
\varepsilon_X
&:=\E\norm{\varphi_1-\varphi_2}_\infty^2\\
&\quad+\E\int_0^T\left(
|b_1(t,(X_2)_t)-b_2(t,(X_2)_t)|^2+
|\sigma_1(t,(X_2)_t)-\sigma_2(t,(X_2)_t)|^2\right)dt,
\\
\varepsilon_Y
&:=\E|\xi_1-\xi_2|^2\\
&\quad+\E\int_0^T
|f_1(t,(X_2)_t,Y_2(t),Z_2(t))
-f_2(t,(X_2)_t,Y_2(t),Z_2(t))|^2dt.
\end{align*}
Assume that $\varepsilon_X$ and $\varepsilon_Y$ are finite. For an Osgood
modulus $\rho$, define
\[
\mathcal B_{\rho,C}(a,\tau)=G_1^{-1}(G_1(a)+C\tau)\quad(a>0),
\qquad \mathcal B_{\rho,C}(0,\tau)=0.
\]
Equivalently, one may extend $G_1$ by $G_1(0)=-\infty$.
Then the differences satisfy
\begin{align}
u(t)&\le C\varepsilon_X+C\int_0^t\rho_X(u(s))ds,\label{eq:stab-X}\\
v(t)&\le C(\varepsilon_X+\varepsilon_Y)
+C\int_t^T\overline\rho_Y(v(s))ds
+C\int_t^T\rho_X(u(s))ds, \label{eq:stab-Y}
\end{align}
where
\[
u(t)=\E\sup_{-\delta\le r\le t}|X_1(r)-X_2(r)|^2,\qquad
v(t)=\E|Y_1(t)-Y_2(t)|^2.
\]
In particular, if $\varepsilon_X+\varepsilon_Y\to0$ along a sequence of
data for which $\rho_X,\rho_Y,L_X,L_f,L_z$ are common to every rank and to
the reference datum, then the constants in the estimates are uniform and
the corresponding solutions converge in
$\Sp\times\Sp\times\Hp$.
More precisely,
\begin{align}
U&:=\mathcal B_{\rho_X,C}(C\varepsilon_X,T),\label{eq:explicit-U}\\
V&:=\mathcal B_{\overline\rho_Y,C}
\bigl(C[\varepsilon_X+\varepsilon_Y+T\rho_X(U)],T\bigr)
\label{eq:explicit-V}
\end{align}
satisfy $\sup_tu(t)\le U$ and $\sup_tv(t)\le V$. Moreover,
\begin{equation}
\begin{split}
&\norm{X_1-X_2}_{\Sp}^2+\norm{Y_1-Y_2}_{\Sp}^2
+\norm{Z_1-Z_2}_{\Hp}^2\\
&\qquad\le U+C\{\varepsilon_X+\varepsilon_Y+T\rho_X(U)
+T\overline\rho_Y(V)\}.
\end{split}\label{eq:full-stability}
\end{equation}
\end{theorem}

\begin{proof}
Subtract the two forward equations. Add and subtract
$b_1(t,(X_2)_t)$ and $\sigma_1(t,(X_2)_t)$. BDG,
Cauchy--Schwarz, the Osgood assumption and Lemma~\ref{lem:delay} give
\eqref{eq:stab-X}. The quantitative Bihari bound gives
\eqref{eq:explicit-U}.

For the backward equations, add and subtract
$f_1(t,(X_2)_t,Y_2,Z_2)$. The proof of
Proposition~\ref{prop:backward-stability}, together with Young's
inequality for the coefficient error, gives \eqref{eq:stab-Y}. The
backward Bihari inequality gives \eqref{eq:explicit-V}. The localized
Itô identity at time zero also yields
\begin{equation}
B:=\E\int_0^T|Z_1(s)-Z_2(s)|^2ds
\le C\{\varepsilon_X+\varepsilon_Y+T\rho_X(U)
+T\overline\rho_Y(V)\}. \label{eq:stab-Z}
\end{equation}
If the perturbations tend to zero, then $U,V,B\to0$, and hence
\[
\sup_{t\le T}\E|Y_1(t)-Y_2(t)|^2+
\E\int_0^T|Z_1(t)-Z_2(t)|^2dt\longrightarrow0.
\]
For completeness, let
\[
A=\E\sup_{0\le t\le T}|Y_1(t)-Y_2(t)|^2,\qquad
M_t=\int_0^t\langle\Delta Y(s),\Delta Z(s)dW(s)\rangle.
\]
Since $2\sup_t|M_T-M_t|\le4\sup_t|M_t|$, BDG and Young's inequality give
\[
4\E\sup_{t\le T}|M_t|\le\frac12A+CB.
\]
Taking the supremum in the localized Itô identity, absorbing $A/2$, and
using the same Osgood estimates yields
\[
A\le C\{\varepsilon_X+\varepsilon_Y+T\rho_X(U)
+T\overline\rho_Y(V)+B\}.
\]
Combining this inequality with \eqref{eq:stab-Z} and the bound for $X$
proves \eqref{eq:full-stability}. Continuity at zero of the two Bihari
comparison maps, for fixed $T$, gives convergence in the full product
space when $\varepsilon_X+\varepsilon_Y\to0$.
\end{proof}

\begin{corollary}[Linear modulus]
If $\rho_X(r)=\rho_Y(r)=r$, Theorem~\ref{thm:main} applies to weakly
monotone drivers satisfying \emph{(O5)}. If, in addition, $f$ is globally
Lipschitz in $y$, then \emph{(O5)} is automatic and one recovers the usual
Lipschitz well-posedness result for the decoupled delayed FBSDE, without a
smallness restriction on $T$.
\end{corollary}

\begin{proof}
The Osgood integral diverges for $\rho(r)=r$, while Bihari's inequality
reduces to Gronwall's lemma. A one-sided Lipschitz inequality alone does
not bound $|f(t,\eta,y,z)-f(t,\eta,0,z)|$; this is why \emph{(O5)} remains
part of the assumptions unless a two-sided global Lipschitz bound is
imposed. Under the uniform bound
$|f(t,0,y,0)-f(t,0,y',0)|\le L_y|y-y'|$, the integrand in
\eqref{eq:local-growth} is bounded by $L_yR$, so \emph{(O5)} is automatic.
\end{proof}

\section{An explicit non-Lipschitz example}

\begin{example}
Let $\alpha$ be the uniform probability measure on $[-\delta,0]$ and
define, in dimension one,
\begin{align*}
b(t,\eta)&=\beta(t)+\vartheta\!\left(
\int_{-\delta}^0\eta(u)\alpha(du)\right),\\
\sigma(t,\eta)&=\gamma(t)+\vartheta\!\left(
\int_{-\delta}^0\eta(u)\alpha(du)\right),\\
f(t,\eta,y,z)&=f_0(t)+a\vartheta\!\left(
\int_{-\delta}^0\eta(u)\alpha(du)\right)+\chi(y)+\ell z,
\end{align*}
where $\beta,\gamma,f_0\in L^2(0,T)$, $a,\ell\in\R$, and
\[
\rho(q)=
\begin{cases}q\log(1/q),&0<q\le e^{-1},\\ e^{-1},&q>e^{-1},\\0,&q=0,
\end{cases}
\qquad \omega(r)=\sqrt{\rho(r^2)},
\]
\[
\vartheta(r)=\omega(r^+),\qquad
\chi(y)=-\operatorname{sgn}(y)\min\{1,\sqrt{|y|}\}.
\]
The function $\omega$ is nondecreasing and concave: on
$(0,e^{-1/2})$ it equals $r\sqrt{-2\log r}$, and it joins the constant
value $e^{-1/2}$ with zero derivative. Indeed, writing
$L=-2\log r$,
\[
\omega'(r)=\sqrt L-L^{-1/2},\qquad
\omega''(r)=-\frac{L^{-1/2}+L^{-3/2}}{r}<0.
\]
Concavity and $\omega(0)=0$ imply
$|\omega(a)-\omega(b)|\le\omega(|a-b|)$ for $a,b\ge0$. Since the positive
part is nonexpansive,
\[
|\vartheta(r)-\vartheta(r')|^2\le\rho(|r-r'|^2).
\]
Moreover, $\chi$ is continuous, bounded and nonincreasing, so
$2(y-y')[\chi(y)-\chi(y')]\le0\le\rho(|y-y'|^2)$.
Jensen's inequality gives
\[
\left|\int_{-\delta}^0(\eta(u)-\eta'(u))\alpha(du)\right|^2
\le\norm{\eta-\eta'}_\alpha^2.
\]
Because $\vartheta(0)=\chi(0)=0$, the free data at the origin are exactly
$\beta,\gamma,f_0$. The boundedness of $\chi$ gives \emph{(O5)} directly.
Thus the coefficients satisfy \emph{(O2)--(O5)} without any additional
truncation, while $\vartheta(r)/r\to\infty$ as $r\downarrow0$. Theorem
\ref{thm:main} gives a unique solution. For instance, the structural
constants may be chosen with $L_X=2$, $L_f\ge a^2$ and
$L_z\ge|\ell|$ (enlarging them to be strictly positive if necessary).
\end{example}

\begin{remark}
The delay in this example concerns the forward segment $X_t$; the
backward variables $Y(t)$ and $Z(t)$ are current values, not delayed
segments.
\end{remark}

\section{Conclusion}

We have established an Osgood well-posedness and stability theory for a
class of FBSDEs with delayed generators. The delayed forward state is
treated on its path space, while the backward driver is one-sided Osgood in
the backward variable and Lipschitz in the martingale integrand. The
forward proof combines delay estimates, Picard iteration and Bihari's
inequality. Backward existence follows from weakly monotone BSDE theory
after the required local integrability is checked explicitly. The joint
stability estimate composes the two Osgood controls.
In particular, the full-norm bound \eqref{eq:full-stability} and the
explicit non-Lipschitz example exhibit the contribution of the triangular
system beyond the two separate well-posedness results.

The forward equation has been kept decoupled in order to obtain a complete
proof under \emph{(O1)--(O5)}, without a smallness restriction on $T$.
The delay in the backward driver acts through $X_t$, whereas $Y(t)$ and
$Z(t)$ are current variables. This differs from models in which histories
of the backward variables enter the generator. Natural continuations are a fully
coupled theorem under an Osgood monotonicity condition, delayed dependence
on the backward history, and $L^p$ estimates for $p\ne2$.

\end{document}